\documentclass[graybox,footinfo]{svmult}

\smartqed
\usepackage{mathptmx}       % selects Times Roman as basic font
\usepackage{helvet}         % selects Helvetica as sans-serif font
\usepackage{courier}        % selects Courier as typewriter font
\usepackage{type1cm}        % activate if the above 3 fonts are
\usepackage{graphicx}       % standard LaTeX graphics tool
\usepackage{array,colortbl}
\usepackage{amsmath,amsfonts,amssymb,bm} % no amsthm, Springer defines Theorem, Lemma, etc themselves

\usepackage{microtype} % good font tricks

\usepackage[colorlinks=true,linkcolor=black,citecolor=black,urlcolor=black]{hyperref}
\usepackage{bookmark}
\makeatletter % to avoid hyperref warnings:
  \providecommand*{\toclevel@author}{999}
  \providecommand*{\toclevel@title}{0}
\makeatother

\usepackage{color}
\usepackage{hyperref} % gives nice possibilities if you check the      
\newcommand{\N}{{\mathbb{N}}} % natural numbers {1, 2, ...}
\newcommand{\R}{{\mathbb{R}}} % reals
\newcommand{\rd}{\,\mathrm{d}} % differential symbol
\newcommand{\norm}[1]{\left\Vert#1\right\Vert}

\usepackage{citeref}   %  generates page numbers, where an item is cited      

\newcommand\novakJAT{{J. Approx. Th.}} 
\newcommand\novakNM{{Numer. Math.}}
\newcommand\novakCA{{Constr.  Approx. }}
\newcommand\novakJC{{J. Complexity }} 
\newcommand\novakMC{{Math.  Comp. }}

\newcommand{\novakchange}[1]{{}{#1}}
\begin{document} 

%%%%%%%%%%%%%%%%%%%%%%%%%%%%%%%%%%%%%%%%%%%%%%%%%%%%%%

% \title*{Numerical Integration and its Complexity: \\ \smallskip
\title*{Optimal Algorithms for Numerical Integration: \\ \smallskip
Recent Results and Open Problems}
% if the original contribution title is too long you may use the following
% short title
\titlerunning{Numerical Integration and its Complexity}
\author{Erich Novak}
\institute{Erich Novak
	\at Mathematisches Institut, University Jena, 
	Ernst-Abbe-Platz 2, D-07743 Jena, Germany\\
 \email{erich.novak@uni-jena.de}}

\maketitle

\abstract{
We present recent 
results on optimal algorithms for numerical integration
and several open problems. 
The paper has \novakchange{six} parts:  
\begin{enumerate} 
    \item \novakchange{Introduction} 
 	\item Lower Bounds  
	\item Universality  
	\item General Domains  
	\item iid Information 
\item Concluding Remarks 
	\end{enumerate}
}

\section{Introduction} 

We study the problem of numerical integration, i.e., of approximating 
the integral
\begin{equation}\label{novakeq01}  
S_d(f) = \int_{D_d} f(x) \rd x
\end{equation}
over an open subset $D_d\subset \R^d$ of 
Lebesgue measure $\lambda^d(D_d)=1$ for integrable functions
$f\colon D_d\to\R$. 
The main interest is on the behavior of the minimal number of function values
that are needed in the worst case setting
to achieve an error of  at most $\varepsilon>0$. 
Note that classical examples of domains $D_d$ are the unit cube $[0,1]^d$ and 
the normalized Euclidean ball (with volume 1),
which are closed. However, we work with their interiors 
for definedness  of certain derivatives. 

We state the problem. Let $F_d$  be a class of %  continuous 
integrable functions $f\colon D_d\to\R$. For~\mbox{$f \in F_d$},
we approximate the integral $S_d(f)$, see~\eqref{novakeq01}, by algorithms
of the form 
$$ 
A_{n}(f)=\phi_{n}(f(x_1),f(x_2),\dots,f(x_n)),  
$$ 
where $x_j \in  D_d$ can be chosen adaptively and $\phi_{n}\colon \R^n\to
\R$ is an arbitrary mapping. Adaption means that the selection of $x_j$
may depend on the already computed values $f(x_1),f(x_2),\dots,f(x_{j-1})$.
We define $N\colon  F_d \to \R^n$ by 
$N(f) = (f(x_1), \dots , f(x_n))$. 
The (worst case) error of the algorithm $A_{n}$ is defined by % \footnote{Some  
% authors still try to assess and compare different algorithms 
% without specifying a class $F_d$ of inputs. This is hopeless  --  in any part 
% of computational mathematics.}
$$
e(A_{n}, F_d )=\sup_{f \in F_d}|S_d(f)-A_{n}(f)|, 
$$
the optimal errors  are given by
$$
e(n, F_d) = \inf_{A_n} \,  e(A_n) .
$$ 
We minimize over all choices of $n$  adaptive sample points $x_j$ and
mappings $\phi_{n}$.
The information complexity $n(\varepsilon,F_d)$ is 
the minimal number of function values 
which is needed to guarantee that the error is
at most $\varepsilon$, i.e., 
$$
n(\varepsilon,F_d)=\min\{n \mid \exists\ A_{n}\ \mbox{such that}\ 
e(A_{n})\le\varepsilon\}.
$$
Here we minimize $n$ over all choices of adaptive sample points $x_j$ and
mappings $\phi_{n}$. 
It is essential to specify a class $F_d$ in order to 
assess and compare different algorithms, 
and to do a complexity analysis. 

It was proved by Smolyak and Bakhvalov 
that as long as the class $F_d$
is convex and symmetric 
we may restrict the minimization of $e(A_n)$ by
considering only nonadaptive 
choices of $x_j$ and  linear mappings $\phi_{n}$, 
i.e.,  it is enough to consider $A_n$ 
of the form
\begin{equation}     \label{novakeq03} 
A_n(f) = \sum_{i=1}^n a_i f(x_i) 
\end{equation}  
and, moreover, 
  \begin{equation} \label{novak1234} 
	e(n, F_d) = %\operatorname*{\vphantom{p}inf}_{x_1,\ldots,x_n}
			\inf_{x_1,\ldots,x_n}
				\sup_{\substack{f \in  F_d \\ N(f) = 0}}
					S_d(f) \, . 
\end{equation} 
\begin{remark} 
a) For a proof of this statement see, 
for example, \cite[Theorem~4.7]{novakNW08}. This result clearly helps  
to find optimal algorithms and to prove  
complexity results. 

b) A~linear algorithm~$A_n$ is called a quasi Monte Carlo (QMC) algorithm 
if~\mbox{$a_i = 1/n$} for all~$i$ and is called a positive 
quadrature formula if~\mbox{$a_i > 0$} for all~$i$. 

c) More on the optimality of linear algorithms 
and on the power of adaption can be found in
\cite{novakCW04,novakNo96,novakNW08,novakTWW88,novakTW80}. 
There are important classes of functions that are \emph{not} 
convex and symmetric,   
see 
\cite{novakCDHHZ,novakGo21,novakKNR19,novakNR96,novakPl15,novakPW09}.
\novakchange{In this paper} we mainly 
consider convex and symmetric $F_d$ 
and then 
we can use formula \eqref{novak1234} for $e(n,F_d)$.

\novakchange{
d) So far we spoke about deterministic algorithms and one 
may ask whether adaption helps for randomized algorithms, 
again under the assumption that $F$ is convex and symmetric. 
This question was open for many years and just recently 
solved by Stefan Heinrich~\cite{novakHe24}:
For \emph{some} classes $F$ adaption helps a lot in the randomized setting.
}
\qed
\end{remark} 

\section{Lower Bounds}

The \emph{optimal rate of convergence} plays an important role 
in numerical analysis.
We start with a classical result of Bakhvalov~\cite{novakBa59} 
for the class 
$$
F^k_d = \{ f\colon  [0,1]^d \to \R \mid 
\Vert D^\alpha f \Vert_\infty \le 1,  \ 
|\alpha| \le k \},
$$
where $k \in  \N$ and $|\alpha| = \sum_{i=1}^d \alpha_i $ for 
$\alpha \in  \N_0^d$ and $D^\alpha f$ denotes the respective 
partial derivative. 
For two sequences $a_n$ and $b_n$ of positive numbers we write 
$a_n \asymp b_n$ if there are positive numbers $c$ and $C$ such 
that 
$c < a_n/b_n < C $ for all $n \in \N$; 
we write $a_n \approx b_n$ if $\lim_{n \to \infty} \frac{a_n}{b_n} = 1$.
Bakhvalov proved that 
\begin{equation}  \label{novak59}  
e(n, F_d^k) \asymp n^{-k/d}  .
\end{equation} 

Observe that we cannot  conclude anything on 
$n(\varepsilon,F_d^k)$ if $\varepsilon$ is fixed and $d$ is large, since 
\eqref{novak59} contains \emph{hidden factors} that depend on 
$k$ and $d$.  Actually the lower bound was of the form
$$
e(n, F_d^k) \ge c_{d,k} n^{-k/d} ,
$$
where the $c_{d,k}$ decrease with $d \to \infty$ and tend to zero. 

%  TECHNIQUE OF BUMP FUNCTIONS. 
%  NEW RESULTS? 

The result of Bakhvalov was recently improved and now we know that 
$n(\varepsilon,F_d^k)$ is of the order $(c_k d/\varepsilon)^{d/k}$.  
A simplified version of the (upper and lower) bounds reads 
as follows. 

\begin{theorem} \cite{novakHNUW17}. \     
\label{novakT1} 
\begin{equation} 
\min(1/2, c_{k,1} d n^{-k/d} ) \le e(n, F_d^k) 
\le \min (1, c_{k,2} d n^{-k/d} )  , 
\end{equation} 
where, for the upper bound, we assume that $n=m^d$, $m \in \N$. 
\end{theorem} 

The lower bounds in \cite{novakHNUW17} are 
explicit and they hold for all domains $D_d$ with volume 1. 
For $k=2$ we have
\begin{equation} 
n(\varepsilon,F_d^2( D_d) ) \ge \frac{2}{2+d} \left( \frac{1}{533}
\frac{d}{d+2} \frac{d}{\varepsilon} 
\right)^{d/2} ,
\end{equation} 
for all $\varepsilon \in (0, d/(d+2))$.
The classes $F_d^k(D_d)$ are defined as the $F_d^k$, the cube $[0,1]^d$ is now replacd 
by a general domain $D_d$. 
 
By the \emph{curse of dimensionality} we mean that 
$n(\varepsilon,F_d)$ is exponentially large in $d$. 
That is, there are positive numbers $c$, $\varepsilon_0$ and $\gamma$ such that
\begin{equation}  \label{novakcurse}    
n(\varepsilon,F_d) \ge c \, (1+\gamma)^d  \quad
\mbox{for all} \quad \varepsilon \le \varepsilon_0  \quad \mbox{and infinitely many} 
\quad  d\in \N. 
\end{equation}

If, on the other hand, 
$n(\varepsilon,F_d)$ is bounded by a polynomial in $d$ and $\varepsilon^{-1}$ 
then we say that the problem is \emph{polynomially tractable}. 
If 
$n(\varepsilon,F_d)$ is bounded by a polynomial in $\varepsilon^{-1}$ alone, i.e., 
$n(\varepsilon,F_d) \le C \varepsilon^{-\alpha}$ for $\varepsilon < 1$
with $C$ and $\alpha$ independent of $d$, 
then we say that the problem is \emph{strongly polynomially tractable}.
 
\begin{remark} 
a) On the proof of Theorem~\ref{novakT1}:
For the upper bound it is 
enough to consider standard product rules that are universal,
they do not depend on the smoothness $k$. 
See Section~\ref{novaks3}  for more on universality. 
New are the lower bounds. 

b) 
It is often said that product rules are bad in higher dimension,
but of course this depends on the function classes. For the classes $F_d^k$ 
product rules are almost optimal. 
   
c) 
It clearly follows from Theorem~\ref{novakT1} 
that the curse of dimensionality holds for the classes $F_d^k$.
It seems that these classes are too big and/or the error criterion 
is too strong. 
\end{remark}

\noindent 
{\bf Open Problem:} \
For which $D_d$ do we have upper bounds of the form
\[
e(n, F_d^k (D_d)) 
\le  c_{k} d n^{-k/d}   ? 
\]
Construct algorithms that achieve this bound for  such 
domains.

\smallskip

We also discuss randomized (or Monte Carlo) algorithms 
and, 
for the purpose of this paper, 
it might be enough that a randomized algorithm $A$
is a random variable $(A^\omega)_{\omega \in \Omega}$ 
with a random element $\omega$ where, 
for each fixed $\omega$, the algorithm $A^\omega$ is a 
(deterministic) algorithm as before.
We denote by $\mu$ the distribution of the $\omega$. 
In addition one needs rather weak measurability assumptions, 
see also the textbook \cite{novakMNR09}. 
Let $\bar n(f,\omega)$ be the number of function values 
used for fixed $\omega$ and $f$. 
The number 
$$
\tilde n(A) = \sup_{f \in F_d} \int_\Omega \bar  n(f, \omega) \rd\mu (\omega)
$$
is called the \emph{cardinality} of the randomized algorithm $A$ and
%  $$
%  e^{\rm ran} (A,F_d) = \sup_{f \in F_d} \left( \int_\Omega^*   
%  \Vert S(f) - \phi_\omega (N_\omega (f))
%  \Vert^2 \rd\mu (\omega) \right)^{1/2}
%  $$
\[
e^{\rm ran} (A,F_d) = \sup_{f \in F_d}  \int_\Omega^*   
\vert S(f) - \phi_\omega (N_\omega (f))
\vert \rd\mu (\omega) 
\]
is the \emph{error} of $A$.
By $\int^*$ we denote the upper integral. For $n \in \N$,  define
$$
e^{\rm ran} (n, F_d) = \inf \{ e^{\rm ran} (A, F_d) \, : \, \tilde n(A) \le n \} .
$$

If $A\colon  F_d \to G$ is a deterministic algorithm
then $A$ can also be treated as a randomized algorithm
with respect to a Dirac (atomic)  measure $\mu$.
In this sense we can say that deterministic algorithms are special
randomized algorithms. Hence the inequality
\begin{equation}
e^{\rm ran} (n, F_d) \le e (n, F_d)
\end{equation}
is trivial. 

\begin{remark} 
The error numbers $e(n,F)$ and $e^{\rm ran} (n,F)$ can be bounded from above 
by the Kolmogorov widths of $F$ in the space 
$L_\infty$ and $L_2$, respectively. 
This is an old result from \cite{novakNo86} that got some attention 
recently since similar upper bounds could be proved 
for several approximation problems 
(optimal recovery in $L_p$ spaces with deterministic or randomized 
algorithms). 
\end{remark}

\begin{remark} 
There is another error criterion for randomized algorithms with two parameters: 
Small error $\varepsilon$ with large probability $1- \delta$. 
The median of mean algorithm is often used, in particular for small $\delta >0$, 
see also 
\cite{novakGL22,novakGS22,novakHJ19,novakKNR19,novakKR19,novakPO22,novakPO22}.  
\end{remark} 

Bakhvalov~\cite{novakBa59} found the optimal rate of convergence
for the numbers $e^{\rm ran} (n, F_d)$  
already in 1959 for the classes 
$
F^k_d = \{ f\colon  [0,1]^d \to \R \mid 
\Vert D^\alpha f \Vert_\infty \le 1,  \ 
|\alpha| \le k \}
$
that we already defined.
He proved
\begin{equation} 
e^{\rm ran} (n, F_d^k) \asymp n^{-k/d-1/2}  .
\end{equation} 

\begin{remark} 
A~proof of the \emph{upper bound} can be given with a technique 
that is often called \emph{separation of the main part} 
or also \emph{control variates}. 
For $n=2m$ use $m$~function values to construct a good $L_2$ approximation 
$f_m$ of $f \in F_d^k$ by a deterministic algorithm. 
The optimal rate of convergence is 
$$
\Vert f - f_m \Vert_2 \asymp m^{-k/d} .
$$
Then use the unbiased estimator
$$
A_n^\omega (f) = S_d (f_m) + \frac{1}{m} \sum_{i=1}^m (f-f_m)(X_i) 
$$ 
with iid random variables $X_i$ that are uniformly distributed 
on $[0,1]^d$. 
See, for example, \cite{novakMNR09,novakNo88} for more details.
An algorithm of Krieg and Novak~\cite{novakKN17} is 
universal, it works for every $k$, see Section~\ref{novaks3}. 
% See also Section~5. 
% 
% We add in passing that the optimal order of convergence can be obtained for 
% many function spaces (Besov spaces, Triebel-Lizorkin spaces) 
% and for arbitrary bounded Lipschitz domains
% $D_d \subset \R^d$; see \cite{novakNT06}, where the approximation problem is 
% studied.
\end{remark} 
 
To obtain an explicit randomized 
algorithm with the optimal rate of convergence 
one needs a random number generator for the set $D_d$. 
If it is not possible to obtain 
efficiently random samples from the uniform distribution 
on $D_d$ one can work with Markov chain Monte 
Carlo (MCMC) methods, see 
Rudolf~\cite{novakRu12}.

\smallskip 
\noindent
{\bf Open Problem: } \
Study, similarly  as in Theorem~\ref{novakT1},  the dependence on the dimension 
for randomized algorithms. 
Using results of Krieg~\cite{novakKr19} one can prove upper bounds 
of the form $d^t \, n^{-k/d-1/2}$. 
One may take $t=(k+1)/2$,  but the optimal (smallest) $t$ is not known.  

\smallskip

Now we study error bounds for quadrature formulas and assume 
that the integrand is from a Hilbert space $F$ of 
real valued functions defined
on a set~$D$. 
We assume that function evaluation is continuous and hence are dealing 
with a reproducing kernel Hilbert space (RKHS) $F$ with a kernel $K$. 
We want to compute $S(f)$ for $f \in F$, where $S$ is a 
continuous
linear functional, 
hence 
$S (f)  =\langle f , h \rangle$ 
for some $h\in F$. 

Vyb{\'\i}ral~\cite{novakVy20}
recently solved old open problems on  
positive definite matrices by improving 
the Schur product theorem. 
This yields new lower bounds for the complexity 
of such integration problems,
see 
\cite{novakHKNV21} \novakchange{and Krieg and  
Vyb{\'\i}ral~\cite{novakKV23}.} 

A main result of 
\cite{novakVy20}
is the following. 
Assume that $M \in \R^{n \times  n}$ is positive semi-definite. 
Then also 
\begin{equation} 
M \circ M - \frac{1}{n} ({\mathop{\mathrm{diag}}} M)({\mathop{\mathrm{diag}}} M)^T
\end{equation} 
is positive semi-definite.
Here $M \circ M$ denotes the pointwise or Schur product 
and 
${\mathop{\mathrm{diag}}} M=(M_{1,1},\dots,M_{n,n})^T$ refers to 
the diagonal entries of $M$ 
whenever $M \in  \R^{n\times n}$. 
%
% Moreover, if $A,B\in\R^{n\times n}$ are two symmetric matrices, 
% we write $A\succeq B$ if $A-B$ is positive semi-definite.
% The Schur product of $A$ and $B$ 
% is the matrix $A\circ B$ with
% \[
% (A\circ B)_{i,j}=A_{i,j}B_{i,j}\quad\text{for}\quad i,j\le n.
% \]
The classical Schur product theorem states that 
the Schur product of two positive semi-definite matrices
is again positive semi-definite.
The result of 
Vyb{\'\i}ral~\cite{novakVy20}
improves this if the two matrices coincide.

Now we use that lower 
bounds for the worst case error of quadrature formulas are equivalent 
to the statement that certain matrices are positive semi-definite:  
The inequality 
\[
e(n, F)^2 \ge \Vert S \Vert^2 - \alpha^{-1}
\]
for an $\alpha >0$ is equivalent to  the fact that 
\[
(K(x_i,x_k) - \alpha h(x_i)h(x_k))_{i,k}
\]
is positive semi-definite for all $x_1, x_2, \dots , x_n$. 

We are mainly interested in 
\emph{tensor product problems}. We will therefore assume that 
$H_i$ is a RKHS
on a domain $D_i$ with kernel $K_i$ for all $i\le d$
and that $F_d$ is the tensor product of these spaces.
That is, $F_d$ is a RKHS on $\mathbf D_d=D_1\times\cdots\times D_d$
with reproducing kernel
\[
 \mathbf K_d: \mathbf D_d\times \mathbf D_d 
 \to \R,\quad \mathbf K_d(x,y)=\prod_{i=1}^d K_i(x^i,y^i).
\] 
If $h_i \in  H_i$ and $\tilde S_i(f)=\langle f,h_i \rangle$ 
for $f\in  H_i$, 
we will denote by $\mathbf h_d$
the tensor product of the functions $h_i$, i.e.,
$$
\mathbf h_d(t)=(h_1\otimes \dots\otimes h_d)(t)=h_1(t^1)\cdot\ldots\cdot h_d(t^d),
\quad t=(t^1,\dots,t^d)\in  \mathbf D_d. %
$$
We study the tensor product functional 
$S_d=\langle\cdot,  \mathbf h_d  \rangle$ on 
the unit ball of  $F_d$.

Before we state a general result, we present the example that initiated 
this research, \novakchange{see \cite{novakNo99}}. 
Here $H_i$ is the three-dimensional RKHS
on $D_i=[0,1]$ with the three orthonormal functions 
$b_1=1$, $b_2=\cos(2\pi x)$ and $b_3=\sin(2\pi x)$,
and we consider the usual integral over $[0,1]^d$.
\novakchange{Vyb{\'\i}ral~\cite{novakVy20} proved,} 
for this example, 
\[
e(n,F_d)^2 \ge 1 - n 2^{-d}
\] 
and positive quadrature formulas are optimal. 
Of course this implies the curse of dimensionality: 
To achieve the modest error $\frac{1}{2} \sqrt{2}$ 
we need $2^{d-1}$ function values. 

A more general result is this: 

\begin{theorem} \cite{novakHKNV21}. \   \label{novakT2} 
For all $i\le d$, let $H_i$ be a RKHS and let $\tilde S_i$ be a bounded linear 
functional on $H_i$ 
with unit norm and nonnegative representer~$h_i$. 
Assume that there are functions $f_i$ and $g_i$ in $H_i$ and a number $\alpha_i\in (0,1]$ such that
$(h_i,f_i,g_i)$ is orthonormal in $H_i$ and $\alpha_i h_i=\sqrt{f_i^2+g_i^2}$.
Then the tensor product problem $S_d=\tilde S_1\otimes\hdots\otimes \tilde S_d$ %
satisfies for all $n\in\N$ that
\[
e(n,F_d)^2 \ge 1-n\prod_{i=1}^d(1+\alpha_i^2)^{-1}.
\]
\end{theorem}

In particular, if all the $\alpha_i$'s are 
equal to some $\alpha>0$ and we want to achieve $e(n, F_d )\le \varepsilon$
for some $0<\varepsilon<1$,
we obtain
\[
n(\varepsilon , F_d )\ge (1-\varepsilon^2)(1+\alpha^2)^d.
\]
This implies the curse of dimensionality.

As an application, we may use this result to obtain lower bounds 
for the complexity of the integration
problem on Korobov spaces with increasing smoothness. 
These lower bounds complement existing upper bounds from~\cite[Section 10.7.4]{novakNW10}.
We add in passing that lower bounds of this form are known and much easier 
to prove for quadrature formulas that only have positive weights, 
see Theorem 10.2 of \cite{novakNW10}.

Using these results 
one can also prove
new lower bounds 
for rather classical Sobolev spaces 
of periodic univariate functions
with small (logarithmic) smoothness, 
see~\cite{novakHKNV22}:
The integration error 
behaves worse than the approximation numbers 
for these spaces.
\novakchange{Multivariate periodic functions and proof techniques 
for lower bounds are studied in the recent paper 
by 
Krieg and  
Vyb{\'\i}ral~\cite{novakKV23}.} 

Lower bounds were  studied  in \cite{novakNW01} with 
the technique of \emph{decomposable kernels}. 
This technique is rather general as long as we consider finite smoothness.
The technique does not work, however, for analytic functions.
In contrast, the new approach 
also works for polynomials and other analytic functions. 
%  One result of  %   \cite{novakHKNV21} 
%  is the following.  
  
\smallskip 
\noindent
{\bf Open Problem: } \ 
The lower bounds of \cite{novakHKNV21} are for Hilbert spaces. 
Do similar lower bounds exist for Banach spaces? 

\section{Universality}  \label{novaks3} 

Let us start with the midpoint rule 
\[
A_n (f) = \frac{1}{n} \sum_{i=1}^n f \left( \frac{2i-1}{2n} \right) 
\] 
for the classes $F^k_1$ defined above. 
The worst case error of $A_n$ is $\frac{1}{4n}$ for $k=1$ and 
the algorithm %E  (more exactly, the sequence $(A_n)_n$ of algorithms) 
is optimal. 
For $k \ge 2$ the worst case error is 
$\frac{1}{24n^2}$, a worst case function is 
$f(x) = \frac{1}{2} (x-\frac{1}{2})^2$. 
For $k=2$ 
the algorithm is still almost optimal, 
we obtain the optimal rate of convergence. 
For $k \ge 3$ the algorithm is bad, the algorithm cannot use 
the high smoothness of $f \in F^k_1$. 
In contrast to this, the Gaussian rules, or the Clenshaw Curtis method, 
are \emph{universal} in the sense that we obtain the optimal 
rate $n^{-k}$ of convergence for any $k$. 

We would like to have more good properties
and consider randomized algorithms that are good 
in the randomized as well as the deterministic 
worst case setting. 

%E  Kann man fuer k=2 Gauss, CC, optimal und midpoint rule vergleichen? 

We start with the simplest case: Periodic functions on $[0,1]$, 
hence $d=1$. 
Then
a very simple algorithm 
has many optimality properties: 
Pick  $n \in  \{ m, m+1, m+2, \dots , 2m \}$ randomly
and then take a randomly shifted midpoint rule 
with $n$ function values and equal weights.
This 
algorithm, but without the random shift, was studied by Bakhvalov~\cite{novakBa61}.

%  Welche Eigenschaften hat dieser Alg.?  Bakhvalov zeigt die 
%  optimale Ordnung fuer Sobolev-Hilbertraeume. Bei Hoelder verliert er was. 

In the case $d=1$ with non-periodic functions 
on $[0,1]$ one can perform a    periodization 
(for example by adding a suitable polynomial) and thus can keep 
the optimal order of convergence but loses  some constants. 
As an alternative, one can take $n=2m$ and first define, via interpolation at the Chebyshev 
points, a polynomial $f_m$ of degree smaller than $m$. Then we use  
separation of the main part together with the simplest Monte Carlo 
method for $(f-f_m)$ to obtain an almost optimal 
algorithm for each class $F^k_1$. 
This algorithm is, asymptotically for large $n$, slightly 
worse (with a logarithmic loss), but possibly  better for 
small $n$. 

\smallskip 
\noindent
{\bf Open Problem: } \
Study and compare these algorithms (for $d=1$) in more detail. 

\smallskip

A randomized algorithm 
\novakchange{$A_n= (A_n^\omega)_{\omega \in \Omega}$} of the form 
\[
A_n^\omega (f) = \sum_{j=1}^n a_j(\omega) f(x_j(\omega))
\]
for $f \in L_1([0,1]^d)$ 
and the integral
$  
S_d(f) = \int_{[0,1]^d} f(x) \rd x
$
was  
introduced  
by Krieg and Novak in \cite{novakKN17}. 
\novakchange{We call this algorithm the ``randomized 
Frolov algorithm''}. 
To state its properties, we define the norms 
\[
\Vert f \Vert^2_{H^k_d} = \sum_{\Vert \alpha \Vert_1 \le k} \Vert 
D^\alpha f \Vert^2_{L_2([0,1]^d)}
\]
and 
\[
\Vert f \Vert^2_{H^{k,{{\rm mix}}}_d } = \sum_{\Vert \alpha \Vert_\infty  \le k} \Vert 
D^\alpha f \Vert^2_{L_2([0,1]^d)} \, ,
\]
see also 
M.~Ullrich~\cite{novakUl17} and  
Krieg~\cite{novakKr19a}.    

\begin{theorem}  \cite{novakKN17,novakUl17}. \     \label{novakT3} 
% There are positive constants $c, c_i$ such that 
% the following holds for all $n \in \N$ 
% with $n \ge c$. 
For all $f \in L_1 ([0,1]^d)$ the 
\novakchange{randomized Frolov} algorithm is \emph{unbiased}, i.e.,  
\begin{equation} 
 {\mathbb{E}} (A_n(f)) = S_d(f) .
\end{equation} 
For classical smoothness $k>d/2$
\[
e^{\rm ran} (A_n, H^{k}_d ) \asymp n^{-k/d-1/2}  \quad \text{and} \quad
e (A_n, H^k_d ) \asymp n^{-k/d} .
\]
For mixed  smoothness $k \in \N$ 
\[
e^{\rm ran} (A_n, H^{k,{{\rm mix}}}_d ) \asymp n^{-k-1/2}  \quad   \text{and}  \quad   
e (A_n, H^{k,{{\rm mix}}}_d ) \asymp n^{-k} (\log n )^{(d-1)/2} .
\]
All these orders of convergence are optimal. 
\end{theorem}

\begin{remark} 
The algorithm from \cite{novakKN17}  has some nice theoretical 
properties\footnote{In particular, the optimal order for 
$e^{\rm ran} (A_n, H^{k,{{\rm mix}}}_d )$, observe that there are no  
log terms, is a new result of 
Mario Ullrich~\cite{novakUl17}.}
but also severe shortcomings:
it is difficult to implement in higher dimensions and
its error is large for moderate values of $n$.
  Hence there is a lot of motivation to 
  look for 
  other approaches to universality.
 
a) The Smolyak algorithm: As long as one only considers deterministic algorithms 
and error bounds
(and the dimension is small), one can take the algorithm from~\cite{novakNR96b}, 
see~\cite{novakNR97}. 
The algorithm is a combination of the Smolyak (sparse grids) algorithm 
with the Clenshaw Curtis quadrature formula. 

For deterministic and stochastic error bounds 
one may use the algorithm of~\cite{novakBNR99} 
(Smolyak algorithm  with interpolation at the Chebyshev points) 
for a  
``separation of the main part'' or control variate and combine it 
with the simplest Monte Carlo method. 
As far as I know, this algorithm was never really propagated 
and studied. 
Compared to the algorithm of \cite{novakKN17}, 
it has slightly worse asymptotic error bounds (with a logarithmic loss), 
but is better for moderate values of $n$ and is easier to implement.
 
Other randomized versions of the Smolyak algorithm were studied 
by Gnewuch and Wnuk~\cite{novakGW20}.
  
b) Lattice rules and  (also randomized) quasi-Monte Carlo methods
based on digital nets
have a great potential for universality. 
We only cite a few papers, together 
with the excellent  books by
Dick and  Pillichshammer~\cite{novakDP08}
and
Dick, Kritzer and Pillichshammer~\cite{novakDKP22}, where the reader 
can find more details:  
\cite{novakDG20},
\cite{novakDGS21},
\cite{novakDP13}, 
\cite{novakEKNO21}, 
\cite{novakGD15},
\cite{novakGS19},
\cite{novakGS22},
\cite{novakGSY16},
\cite{novakKKNU19},
\novakchange{\cite{novakKNW23}}.
See also other papers by Dick, Goda, Kritzer, Kuo, Nuyens, Pillichshammer, Suzuki  
and their coworkers. 
\end{remark} 

\noindent
{\bf Open Problem: } \ 
Of course one should study universality for 
tractability questions, not only for the rate of convergence.  
Some of the cited papers already contain results in this direction.
See, in particular, Chapter~11 of Dick, Kritzer and Pillichshammer~\cite{novakDKP22} 
and 
the papers 
\cite{novakDGS21,novakGS22,novakKKNU19,novakKNW23}. 

\smallskip 
\noindent
{\bf Open Problem: } \ 
So far only the domain $[0,1]^d$ was studied for 
the construction of universal algorithms.  
One might speculate that universal algorithms also 
exist for other domains. 

\section{General Domains}

We usually assumed that the domain is a cube, $D_d=[0,1]^d$, 
now we discuss results for 
more general domains 
$D_d$.
We study  different concepts to measure the complexity.

\subsection*{Order of convergence} 

It is known for many problems that the order of convergence 
of optimal algorithms does not depend 
on the domain $D_d \subset \R^d$. 
Results for Lipschitz domains (and beyond) 
and classical (isotropic) function spaces 
can be found in 
\cite{novakKNS21,novakNWW04,novakNT06,novakWSH22}. 
The paper by Wenzel, Santin and Haasdonk~\cite{novakWSH22} 
has results saying that the order of convergence cannot decrease 
if $D_d$ is replaced by a smaller set. 
Such results depend on the way how we define 
Sobolev (and related) spaces: 
Usually all definitions coincide if $D_d$ is a bounded Lipschitz 
domain but they do not coincide for more general domains. 

Function recovery and integration on manifolds are studied  by 
Brandolini, Choirat, Colzani, Gigante, Seri 
and Travaglini~\cite{novakBCCGST14},
Ehler, Graef and Oates~\cite{novakEGO19}
and 
Krieg and Sonnleitner~\cite{novakKS21}. 
For spaces of functions 
with dominating mixed smoothness many 
problems are still open, see Triebel~\cite{novakTr19}. 

\subsection*{Asymptotic constant and explicit uniform bounds} 

If the rate of convergence does not depend on the domain 
$D_d$, then we may ask about the asymptotic constants
and we may study 
explicit uniform bounds, i.e., 
bounds that hold for all $n$ (number of pieces of information) 
and all (normalized) domains. 

We 
present examples for which the following statements 
are true: 
\begin{enumerate} 
\item 
Also the asymptotic constant does not depend on the shape of $D_d$
or the imposed boundary values,  
but only  on the volume of the domain.
These results can be compared with results of 
Weyl~\cite{novakWeyl} who 
proved that the asymptotic constant for the size of the eigenvalues 
of the Dirichlet Laplacian 
and also of the Neumann Laplacian do not depend on the shape of $D_d$,
it only depends on the volume of $D_d$.
\item
There are explicit and uniform lower (or upper, respectively) 
bounds for the error that are 
only slightly smaller (or larger, respectively) 
than the asymptotic error bound. 
%  By explicit we mean that the bounds hold for any $n$, not just asymptotically. 
%  By uniform we mean that they hold for all domains of a given size. 
\end{enumerate} 

General bounded Lipschitz domains are studied in 
\cite{novakNWW04,novakNT06,novakTr08,novakTr10,novakVy06,novakVy07b}.
In many cases, the optimal order of the errors  is of the form 
$$
e (n, F_d)  \asymp n^{-\alpha} \, (\log n)^\beta ,  
$$
where $\alpha$ and $\beta$ do not depend on the domain $D_d$. 
It is interesting to know also the exact asymptotic constant 
$$
C: = \lim_{n\to \infty}  e(n, F_d)   \, n^{\alpha} \, (\log n)^{-\beta} ,
$$
if it exists. 
The value of $C$ is known only in rare cases 
(unless $d=1$, we do not discuss the univariate case 
in detail), 
and usually only for very special domains, like the cube, see 
\cite{novakCKS16,novakKr18,novakKSU14,novakKSU15}. 

We assume that $D_d \subset \R^d$ is bounded and a metric on $D_d$ is induced 
by a norm in $\R^d$. We denote by $B$ the unit 
ball (in $\R^d$)  and also write $\Vert \cdot \Vert_B$ for the norm. 
To simplify the formulas we consider 
only Lipschitz functions. 
We denote the respective \novakchange{set} by $F^B(D_d)$; 
it contains all functions $f\colon   D_d \to \R$  with 
$$
|f(x) - f(y)| \le \Vert x - y \Vert_B.
$$

We study the problem of numerical integration
and assume that $D_d \subset \R^d$ is Jordan measurable with 
$0 < \lambda^d (D_d) < \infty$. 
Asymptotic formulas, even for more general (weighted)  integration problems, 
were proved by Chernaya~\cite{novakCh95} and by Gruber~\cite{novakGr04}. 
We add an upper bound for the asymptotic constant $\xi_B$  and see that 
it is very close to the lower bound, in particular for large $d$. 
The explicit uniform lower bound is from \cite{novakNo20}. 

\begin{theorem}  
Assume that $D_d \subset \R^d$ is a bounded Jordan measurable 
set with an interior point. 
	Then 
	\begin{equation} 
		e(n, F^B(D_d) )   \approx \xi_B  \,  \lambda^d(D_d)   
		\left( \frac{\lambda^d(D_d)}{\lambda^d(B)}\right)^{1/d} 
	\cdot n^{-1/d},
	\end{equation} 
	where  $\xi_B$  is a constant that depends on the norm and 
\[
	\frac{d}{d+1}  \le \xi_B \le     %   \frac{d}{d+1}  \Theta_B^{1/d}  \le    
	 \frac{d}{d+1} ( d \log d + d \log \log d + 5d )^{1/d}  .
\]  
	Moreover, 
	\begin{equation}
\inf_{D_d}   e(n, F^B(D_d)  ) \cdot \lambda^d  (D_d)^{-(d+1)/d} = \frac{d}{d+1}  \lambda^d (B)^{-1/d} 
		\cdot n^{-1/d} ,  
	\end{equation} 	
the infimum is attained when $D_d$ is  the disjoint union 
of $n$ balls with equal radius.
\end{theorem} 

A similar result for Sobolev spaces
$W^1_p (D_d) $ with $p>d$  was recently obtained by 
 Babenko, Babenko and Kovalenko~\cite{novakBBK21}.

We now discuss  $C^k$-functions.
%   We start with a result from \cite{HNUW17}.
We know from \cite{novakHNUW17}     
the following statement
for the classes
$
F^k_d (D_d)  = \{ f\colon  D_d  \to \R \mid 
\Vert D^\alpha f \Vert_\infty \le 1,  \ 
|\alpha| \le k \}
$ for all 
domains $D_d$ with volume 1:
For all $k \in \N$ there exists a constant $c_k>0$ such that for all 
$d,n \in \N$ 
$$
e ( n,  F^k_d  (D_d) ) \ge \min \{ 1/2, \,  c_k \, d \, n^{-k/d} \} . 
$$ 

%  Observe that the constant $c_k$ does not depend on $D_d$ or $d$, we have a 
%  lower bound for  all $D_d$ with volume 1. 

\begin{remark} 
1) The lower bound cannot be improved since for cubes we have a similar upper bound. 
It would be good to know more on the constants $c_k$ and on extremal sets $D_d$,
where $e( n,  F^k_d(D_d)  )$ is small  for given $k$, $d$ and $n$. 

For these function spaces there cannot be a meaningful explicit uniform  upper bound 
since
$\sup_{D_d}  e( n,  F^k_d(D_d) ) =1$. 
The supremum over $D_d$ makes sense if we impose boundary conditions 
such as $f(x) = 0$ for $x \in \partial D_d $.  

2) 
It is known that 
the weak order
$$
e( n,  F^k_d(D_d)  ) \asymp n^{-k/d} 
$$
holds at least for every bounded Lipschitz domain. 
This follows from  more general results of \cite{novakKNS21,novakNWW04,novakNT06}.  
We guess that more is true and the asymptotic constant 
$$
\lim_{n \to \infty} e(  n,  F^k_d(D_d)  )  \cdot \novakchange{n^{k/d}} = C_{D_d} 
$$
does not depend on $D_d$ (for fixed $d$) if $D_d$ is Jordan measurable
with $\lambda^d(D_d)=1$; 
possibly under some assumptions on the boundary of $D_d$. 
\end{remark} 

\begin{remark}
The norm 
$$
\max_{|\beta|_1 \le k} \Vert D^\beta f \Vert_\infty 
$$
might be reasonable when we consider a cube $D_d = [0,1]^d$
but it is not invariant with respect to rotation. 
The function $f(x) = \sum x_i$ has a gradient with length 
$d^{1/2}$ as $g(x) = d^{1/2} \cdot x_1$, but all 
partial derivatives of $f$ are bounded by one. 
Therefore we also consider an orthogonal invariant norm. 
Since for this modified space $\tilde C^k(D_d)$ 
many problems are still open, we only discuss the case $k=2$ 
in the following.
\end{remark} 

\begin{example} 
Let us discuss the integration problem for the class 
$$
\tilde  C^2 (D_d) = \{ f\in C^2 (D_d) \mid 
\Vert f \Vert_\infty \le 1, \, 
{{\rm Lip}} (f) \le d^{-1/2}, \, 
{{\rm Lip}} (D^\Theta f ) \le d^{-1} \}, 
$$
as in \cite{novakHNUW12,novakHNUW14b,novakHPU19}.
Here $D^\Theta f$ denotes any directional derivative in a direction 
$\Theta \in S^{d-1}$ and 
$$
{{\rm Lip}} (g) = \sup_{x,y\in D_d} \frac{|g(x)-g(y)|}{\Vert x-y\Vert_2} .
$$ 
	Again we assume that $D_d \subset \R^d$ is an open set 
	with $\lambda^d(D_d) = 1$. 
We conjecture that there exists a constant $C>0$ (independent on $d$ and $D_d$) such that 
\begin{equation}  \label{novak123} 
e (n,  \tilde   C^2 (D_d) ) \ge C \cdot n^{-2/d} .
\end{equation} 
\end{example} 

Some comments are in order: 

1) This would be another explicit uniform lower bound and, because of known upper bounds 
for the cube, it certainly cannot be improved. 

2) It is known that the integration problem for $\tilde C^2 (D_d)$ and 
\emph{certain} $D_d$ suffers from the curse 
of dimensionality. This is true if $D_d$ has a small radius, see \cite{novakHPU19} for the 
best known results. 
For example, the curse is known if $D_d$ is a $\ell_p^d$ ball and $p \ge 2$. 
It is not known for $p$ balls and $p<2$.  

3) It is easy to see that the lower bound \eqref{novak123} is true if 
$D_d$ is a disjoint union of $n$  Euclidean balls of the same size. 
But, of course, this domain is not extremal for the given norm. 

4) Assume that the $x_i$ form a grid
and $D_d$ is a cube. Then, for $d=1$, one may take a 
quadratic spline $f_1$ as a fooling function and for $d>1$ one can take a 
fooling function of the form 
$$
f_d(x) = \frac{1}{d} \sum_{i=1}^d f_1 (x^i) ,
$$
where $x= (x^1, \dots , x^d)$. 
	Hence we obtain exactly the lower bound \eqref{novak123} in this case. 
	Loosely speaking, 
	the conjecture says that arbitrary $D_d$ and function values 
	at arbitrary points $x_1, x_2, \dots , x_n$  do not lead to much better 
	error bounds than grids for cubes. 

\section{iid Information} 

So far we mainly studied (deterministic or randomized) 
quadrature formulas
\[
A_n(f) = \sum_{i=1}^n a_i f(x_i) 
\]
that are, in some way, optimal. 
The underlying assumption is that we, when we design an algorithm, 
can choose the points $x_i$ and
the weights $a_i$. 
This is a standard assumption in numerical analysis and complexity 
theory, but it is not always realistic. 

In some parts of  learning,   data science and uncertainty quantification
it is more realistic to assume that the $x_i$ are random and realizations of 
iid random variables, 
see~\cite{novakBGKP21,novakSB14,novakSC08,novakZh20}.
Then we cannot influence the data stream
\[
(x_1, f(x_1), x_2, f(x_2), x_3, f(x_3), \dots ) 
\] 
and, of course, the construction of ``optimal'' 
points $x_1, x_2, \dots , x_n$ 
does not make sense;  
we simply have to use 
the 
information as it comes in. 

A natural question to ask is 
``What is the power of random information?'',
with two possible answers: 
Depending on the assumptions and involved 
function spaces, 
random information might be almost as good as optimal 
information -- or much worse.
In the framework of information-based complexity,
the power of iid information was recently studied and surveyed 
in~\cite{novakHKNPUsurvey,novakHKNPU21a,novakHJ19,novakKNR19,novakSo22}.

I summarize 
results of the paper  
Krieg, Novak and Sonnleitner~\cite{novakKNS21}
where it is proved, formulated as a slogan, that 
random information is quite often 
(almost) as good as optimal information. 

We study integration for functions 
from the Sobolev space $W^s_p(D_d)$ and compare optimal 
randomized (Monte Carlo) algorithms with algorithms 
that can only use iid sample points, uniformly distributed on the domain. 
The main result is that we obtain the same optimal rate 
of convergence if we restrict ourselves  to iid sampling. 

Let $D_d  \subset {\mathbb{R} }^d$ be open and bounded.
We assume that $D_d $ satisfies an interior cone condition.
We also  study the problem of approximating
a function $f$ from the Sobolev space $W_p^s(D_d )$
in the $L_q(D_d )$-norm based on function values $f(x_j)$ 
on a finite set of sampling points $P=\{x_1,\hdots,x_n\}$.
This makes sense if $s > d/p$, 
in which case $W_p^s(D_d)$ is compactly embedded
into the space of bounded continuous functions~$C_b(D_d )$, 
and in the case $p=1$ and $d=s$.  
There is a vast literature on the error for optimal sampling points. 
It is known that the rate of convergence of the 
worst case error of optimal 
deterministic algorithms is
$
 n^{-s/d}
$
for the integration problem
and
$
 n^{-s/d+(1/p-1/q)_+}
$ 
for the approximation problem. 
For general domains, we refer to Narcowich, Wendland and 
Ward \cite{novakNWW04} as well as \cite{novakNT06}
and \cite{novakKNS21}. 

Now  we assume that  the sampling points
are given to us as realizations 
of independent random variables
which are uniformly distributed on the domain.
Then  $P \subset D_d $ is a realization of a \novakchange{distribution} $\mathcal{P}$ and  
also an algorithm is 
a randomized or Monte Carlo algorithm and one may compare 
the algorithm with optimal randomized algorithms
based on 
the expected  error
for inputs from the unit ball of $W_p^s(D_d )$. 

The class of all randomized algorithms is very large 
since one may construct the random variable $x_{k+1}$ 
based on the realization 
$(x_1, f(x_1), \dots , x_k, f(x_k))$
in a sophisticated and complicated adaptive way. 
Nevertheless, 
using this weaker notion of error, 
one cannot improve the rate  
$
 n^{-s/d+(1/p-1/q)_+}
$ 
for the approximation problem, 
see Math\'e~\cite{novakMa91} and Heinrich~\cite{novakHe08} 
as well as Heinrich~\cite{novakHe09a,novakHe09b} for other Sobolev embeddings.
For the integration problem one now obtains the improved 
order
$n^{-s/d-1/2}$ if $p \ge 2$ 
and 
$n^{-s/d-1+1/p}$ if $1 \le p <2$, 
see Bakhvalov~\cite{novakBa62} and \cite{novakNo88}.
 
With algorithms 
based on iid  sampling 
we still obtain the optimal order 
of convergence for the approximation problem 
unless $p=q=\infty$: In this limiting case 
there is a logarithmic loss. 
For the integration problem we obtain the optimal order 
for all $p$. 

We describe the  results in more detail. 
We assume that $D_d \subset {\mathbb{R} }^d$ is open and bounded
and satisfies an interior cone condition.
That is, there is some $r>0$ and $\theta \in (0,\pi]$ 
such that for every $x\in D_d $, we find a cone 
\[
K(x):=\{x+\lambda y: y\in\mathbb{S}^{d-1}, \langle y,\xi(x)\rangle\ge\cos\theta,\lambda
\in  [0,r] \}
\]
with apex $x$, direction $\xi(x)$ in the unit sphere $\mathbb{S}^{d-1}$, 
opening angle $\theta$ and radius $r$
such that $K(x) \subset D_d $. 
Here  $\langle \cdot, \cdot \rangle$ denotes the standard inner product.
For $s\in {\mathbb{N} }$ and $1\le p \le \infty$
such that $s > d/p$ or $p=1$ and $s=d$, we consider the Sobolev space
\[
W_p^s(D_d ) \,:=\, \left\{ f \colon D_d  \to {\mathbb{R} }  \ \big\vert\ D^\alpha f 
 \in L_p(D_d ) \text{ for all } \alpha \in {\mathbb{N} }_0^d \text{ with } |\alpha| \le s  \right\}
\]
with semi-norm
\[
 |f|_{W_p^s(D_d )} \,:=\, \left(\sum_{|\alpha| = s} \Vert D^\alpha f \Vert_{L_p(D_d )}^p\right)^{1/p}
\]
and norm 
\[
 \Vert f\Vert_{W_p^s(D_d )} := 
 \left(\sum_{|\alpha| \le s} \Vert D^\alpha f \Vert_{L_p(D_d )}^p\right)^{1/p},
\]
with the usual modification for $p= \infty$.
Note that $W_p^s(D_d )$ is continuously embedded into $C_b(D_d )$, 
the space of bounded continuous functions 
with the sup norm, and that the embedding is compact
in the case $s>d/p$, see e.g.\ Maz'ya \cite[Section~1.4]{novakMaz11} 
for this fact for domains satisfying an interior cone condition.
We study the problem of $L_q(D_d )$-approximation ($1\le q \le \infty$)
on $W_p^s(D_d )$. 

For a random operator $A\colon W_p^s(D_d ) \to L_q(D_d )$
we study the Monte Carlo error
\[
e^{\rm ran}     \left(A,W_p^s(D_d ),L_q(D_d )\right)
\,:=\,
\sup_{\Vert f\Vert_{W_p^s(D_d )} \le 1}\, {\mathbb{E}}\ \left\Vert f - A(f) \right\Vert_{L_q(D_d )}.
\]
We 
may use the notation ${\mathbb{E}}$ for any mapping if we take the upper integral
for the definition.  
Given a distribution $\mathcal{P}$, we put
\[
e^{\rm ran}    \left(\mathcal{P},W_p^s(D_d ),L_q(D_d )\right)
\,:=\,
\inf_A\, e\left(A,W_p^s(D_d ),L_q(D_d )\right),
\]
where the infimum is taken over all random operators 
of the form $A(f)=\varphi(f|_P)$ 
with a random mapping $\varphi\colon {\mathbb{R} }^P \to L_q(D_d)$. 
Note that $f|_P=\big(f(x)\big)_{x\in P}$ is the restriction 
of $f$ to the point set $P$ and we use this as our information. 
Here $P$ is a realization of $\mathcal{P}$.
This is the smallest Monte Carlo error that can be achieved with the
\novakchange{distribution}  $\mathcal{P}$.
Moreover, we put
\[
e^{\rm ran}    \left(n,W_p^s(D_d ),L_q(D_d )\right)
\,:=\,
\inf_{\mathcal{P}}\, e\left(\mathcal{P},W_p^s(D_d ),L_q(D_d )\right),
\]
where the infimum is taken over all random point sets
of cardinality at most~$n$.
This is the smallest Monte Carlo error that can be achieved with
$n$ optimally chosen random  sampling points.
It is known, at least for special domains $D_d $,  that
\begin{equation}
 e^{\rm ran}    \left(n,W_p^s(D_d ),L_q(D_d )\right) \,\asymp\, n^{-s/d+(1/p-1/q)_+},
\end{equation}
see Math\'e~\cite{novakMa91} and Remark~\ref{novakrem1}.

\begin{theorem}\label{novakthm:main}\cite{novakKNS21}. \ 
 Let $D_d  \subset {\mathbb{R} }^d$ be open and bounded, satisfying an interior cone condition
 and let $1 \le p,q \le \infty$ and $s\in {\mathbb{N} }$ such that $s>d/p$ or $p=1$ and $s=d$.
 For every $n\in {\mathbb{N} }$, let $\mathcal{P}_n$ be a set of $n$ independent and uniformly distributed points 
 on~$D_d $.
 Then
 \[
 e^{\rm ran}    \left(\mathcal{P}_n,W_p^s(D_d ),L_q(D_d )\right)
 \,\asymp\,
 \begin{cases}
\, \displaystyle \left(n / \log n \right)^{-s/d} 
& \text{if } p=q=\infty,\\
\, n^{-s/d+(1/p-1/q)_+} 
& \text{else}. \vphantom{\Big|}
\end{cases}
\]
\end{theorem}
This means that independent and uniformly distributed points 
are (asymptotically) as good as optimally selected 
(deterministic or random) sampling points
in all cases except $p=q=\infty$.

This answers the question for the power of independent uniformly distributed samples
with respect to the Monte Carlo error criterion.
On the other hand, one might also be interested 
in a stronger uniform error criterion.
For a random operator $A\colon W^{s}_{p}(D_d)\to L_{q}(D_d)$ 
the uniform error can be defined by
\[
	{e^{\rm unif}}\left(A,W^{s}_{p}(D_d),L_{q}(D_d)\right)
	\,=\, {\mathbb{E}} \sup_{\norm{f}_{W^{s}_{p}(D_d)}\le 1}\norm{f-A(f)}_{L_{q}(D_d)}.
\]
Note that the order of the supremum and the expected value is interchanged.
Thus, while a small Monte Carlo error \novakchange{$e^{\rm ran}$} means
that for every individual function
the error is small with high probability, 
a small uniform error ${e^{\rm unif}}$ means 
that with high probability
the error is small for every function. 
As before, given a \novakchange{distribution}  $\mathcal{P}$,  we put
\[
{e^{\rm unif}}\left(\mathcal{P},W_p^s(D_d ),L_q(D_d )\right)
\,:=\,
\inf_A\, {e^{\rm unif}}\left(A,W_p^s(D_d ),L_q(D_d )\right),
\]
where the infimum is taken over all random operators 
of the form  $A(f)=\varphi(f|_P)$ 
with a random mapping $\varphi\colon {\mathbb{R} }^P \to L_q(D_d )$. 
Again $P$ is a realization of $\mathcal{P}$. 
This is thus the smallest uniform error that can be achieved
with the \novakchange{distribution} $\mathcal{P}$.
Moreover, we put
\[
{e^{\rm unif}}\left(n,W_p^s(D_d ),L_q(D_d )\right)
\,:=\,
\inf_\mathcal{P}\, {e^{\rm unif}}\left(\mathcal{P},W_p^s(D_d ),L_q(D_d )\right),
\]
where the infimum is taken over all random point sets
of cardinality at most $n$. 
It is known that 
\begin{equation}
 {e^{\rm unif}}\left(n,W_p^s(D_d),L_q(D_d )\right)
 \,\asymp\, n^{-s/d+(1/p-1/q)_+},
\end{equation}
see again \cite{novakKNS21,novakNWW04,novakNT06}.

For the uniform error of independent uniformly distributed samples,
the following result has been obtained in Krieg and Sonnleitner~\cite{novakKS20} 
for bounded convex domains.
See also Ehler, Graef and Oates~\cite{novakEGO19}
and the survey \cite{novakHKNPUsurvey} for earlier results in this direction.

\begin{theorem}\label{novakthm:main-det}\cite{novakKNS21}. \ 
 Let $D_d \subset {\mathbb{R} }^d$ be open and bounded, satisfying an interior cone condition
 and let $1 \le p,q \le \infty$ and $s\in {\mathbb{N} }$ such that $s>d/p$ or $p=1$ and $s=d$.
 For every $n\in {\mathbb{N} }$, let $\mathcal{P}_n$ be a set of $n$ independent and uniformly 
 distributed points on~$D_d$.
 Then
 \[
  {e^{\rm unif}}\left(\mathcal{P}_n,W_p^s(D_d ),L_q(D_d )\right)
 \,\asymp\, 
 \begin{cases}
\, \displaystyle \left( n / \log n \right)^{-s/d+1/p-1/q}  
& \text{if } q\ge p,\\
\ n^{-s/d}
& \text{if } q<p. \vphantom{\Big|}
\end{cases}
\]
\end{theorem}
This means, with respect to the uniform error criterion,
independent and uniformly distributed samples are 
as good as optimally chosen sampling points if and only if $q<p$.
Moreover, if we are bound to independent and uniformly distributed samples,
one can achieve the same rate for the uniform error as for the Monte Carlo error
if and only if $q<p$ or $p=q=\infty$.
In all other cases,
the Monte Carlo error criterion provides a speed-up in comparison
to the uniform error criterion.

All the upper bounds of Theorems~\ref{novakthm:main} and \ref{novakthm:main-det} 
are achieved by the same algorithm which is based on a polynomial 
reproducing map.
The algorithm works for all $p$ and $q$ and up to a given 
smoothness $s$,  
but should be simplified for an implementation.

Let us now turn to the integration problem
\[
S_d (f) = \int_{D_d } f(x) \, {\rm d}  x\,.
\]
We use a similar notation.
For a random operator $A\colon W_p^s(D_d ) \to {\mathbb{R} }$
we define the Monte Carlo error
\[
e^{\rm ran}   \left(A,W_p^s(D_d ), S_d \right)
\,:=\,
\sup_{\Vert f\Vert_{W_p^s(D_d )} \le 1}\, {\mathbb{E}}\ \left\vert S_d (f) - A(f) \right\vert.
\]
Given a \novakchange{distribution} $\mathcal{P}$, we put
\[
e^{\rm ran}    \left(\mathcal{P},W_p^s(D_d ), S_d  \right)
\,:=\,
\inf_A\, e\left(A,W_p^s(D_d ), S_d  \right),
\]
where the infimum is taken over all random operators 
of the form $A(f)=\varphi(f|_P)$ 
with a random mapping $\varphi\colon {\mathbb{R} }^P \to {\mathbb{R} } $. Again  $P$ is a realization of $\mathcal{P}$. 

This is the smallest Monte Carlo error that can be achieved with the
\novakchange{distribution}   $\mathcal{P}$.
Moreover, we put
\[
e^{\rm ran}  \left(n,W_p^s(D_d ), S_d  \right)
\,:=\,
\inf_\mathcal{P}\, e\left(\mathcal{P},W_p^s(D_d ), S_d  \right),
\]
where the infimum is taken over all random point sets
of cardinality at most $n$.
This is the smallest Monte Carlo error that can be achieved with
$n$ optimally chosen sampling points.
It is known, at least for special domains $D_d $,  that 
\[
 e^{\rm ran}   \left(n,W_p^s(D_d ), S_d \right) \,\asymp\, n^{-s/d+ (1/p-1/2)_+ - 1/2 }, 
\]
see Bakhvalov~\cite{novakBa59,novakBa62} and \cite{novakNo88}. 
Again,  we are interested in the smallest Monte Carlo error
which can be achieved with $n$ independent and uniformly distributed
sampling points. As a corollary to Theorem~\ref{novakthm:main} we obtain the following.

\begin{corollary}   \cite{novakKNS21}. \ 
 Let $D_d  \subset {\mathbb{R} }^d$ be open and bounded, satisfying an interior cone condition
 and let $1 \le p \le \infty$ and $s\in {\mathbb{N} }$ such that $s>d/p$
 or $p=1$ and $s=d$.
 For every $n\in {\mathbb{N} }$, let $\mathcal{P}_n$ be a set of $n$ independent and uniformly distributed points 
 on~$D_d $.
 Then
 \[
 e^{\rm ran}   \left(\mathcal{P}_n,W_p^s(D_d ), S_d \right)
 \,\asymp\,
 n^{-s/d+ (1/p-1/2)_+ -1/2 } . 
\]
\end{corollary}

This means that, with respect to the Monte Carlo error, 
independent and uniformly distributed points 
are (asymptotically) as good as optimally selected 
(deterministic or random) sampling points
in all cases. 

This answers the question for the power of independent uniformly distributed samples
for numerical integration
with respect to the Monte Carlo error criterion.
Again, one might also be interested 
in a stronger uniform error criterion.
For a random operator $A\colon W^{s}_{p}(D_d )\to  {\mathbb{R} } $ 
the uniform error can be defined by
\[
	{e^{\rm unif}}\left(A,W^{s}_{p}(D_d ),  S_d  \right)
	\,=\, {\mathbb{E}} \sup_{\norm{f}_{W^{s}_{p}(D_d )}\le 1}  \vert S_d (f) - A(f) \vert .
\]
Again the order of the supremum and the expected value is interchanged.
As before, given a \novakchange{distribution}  $\mathcal{P}\subset D_d $, we put
\[
{e^{\rm unif}}\left(\mathcal{P},W_p^s(D_d ), S_d \right)
\,:=\,
\inf_A\, {e^{\rm unif}}\left(A,W_p^s(D_d ), S_d  \right),
\]
where the infimum is taken over all random operators 
of the form $A(f)=\varphi(f|_P)$ 
with a random mapping $\varphi\colon {\mathbb{R} }^P \to {\mathbb{R} } $. Here  $P$ is a realization of $\mathcal{P}$. 

This is thus the smallest uniform error that can be achieved
with the \novakchange{distribution} $\mathcal{P}$.
Moreover, we put
\[
{e^{\rm unif}}\left(n,W_p^s(D_d ),  S_d  \right)
\,:=\,
\inf_\mathcal{P}\, {e^{\rm unif}}\left(\mathcal{P},W_p^s(D_d ),  S_d  \right),
\]
where the infimum is taken over all random point sets
of cardinality at most $n$. 

For the uniform error of independent uniformly distributed samples,
the following result has been obtained in Krieg and Sonnleitner~\cite{novakKS20} 
for bounded convex domains.

\begin{corollary}   \cite{novakKNS21}. \ 
 Let $D_d  \subset {\mathbb{R} }^d$ be open and bounded, satisfying an interior cone condition
 and let $1 \le p \le \infty$ and $s\in {\mathbb{N} }$ such that $s>d/p$
 or $p=1$ and $s=d$.
 For every $n\in {\mathbb{N} }$, let $\mathcal{P}_n$ be a set of $n$ independent and uniformly distributed points
 on~$D_d $.
 Then
 \[
  {e^{\rm unif}}\left(\mathcal{P}_n,W_p^s(D_d ), S_d \right)
 \,\asymp\, 
 \begin{cases}
\, \displaystyle \left( n / \log n \right)^{-s/d}  
& \text{if } p=1 ,\\
\, n^{-s/d}
& \text{if } p>1 . \vphantom{\bigg|}
\end{cases}
\]
\end{corollary}

This means, with respect to the uniform error criterion,
independent and uniformly distributed samples are 
as good as optimally chosen sampling points 
if and only if $p>1$.
We also notice that for independent and uniformly distributed samples,
the Monte Carlo error criterion provides a speed-up in comparison
to the uniform error criterion for all $p$. 

\begin{remark}  \label{novakrem1} 
Different authors study different domains and possibly even 
different function spaces $W^s_p(D_d )$. 
All the standard definitions of the Sobolev spaces 
coincide if $D_d $ is a bounded Lipschitz domain;
if the domain is not Lipschitz then different texts 
possibly use different spaces. 
If the Sobolev space on $D_d $ is defined by restriction 
of the functions from $W^s_p({\mathbb{R} }^d)$ then one 
obtains smaller spaces than the spaces defined above. 
The stated   results stay true for this altered definition 
since the lower bounds
already hold for functions with compact support inside~$D_d $.
\end{remark}

\begin{remark}
%  Our techniques 
%  can most likely 
%  be used to prove similar results
%  for more general function spaces of isotropic smoothness 
%  like Triebel-Lizorkin or Besov spaces as well as Sobolev spaces on 
%  manifolds.
Another interesting family of spaces are function spaces
of mixed smoothness as surveyed in
D\~ung, Temlyakov and T.~Ullrich~\cite{novakDTU16}.
Here  we are still  far from understanding
the power of 
independent and uniformly distributed sampling points,
and even the power of optimal sampling points
is not known in all  cases.
There are recent results in this direction for 
the special case of $L_2$-approximation 
on the Hilbert space
$ H^{k,{{\rm mix}}}(\mathbb T^d)$
of functions with mixed smoothness $k$ on the $d$-torus.
Namely, it is known that independent and uniformly distributed sampling points 
are optimal with respect to the Monte Carlo error~\cite{novakKri19}
and optimal up to a logarithmic factor with respect to the
uniform error~\cite{novakDKU22,novakKU19,novakUll20}.
See also 
\cite{novakBSU22,novakCD21,novakKU21,novakNSU20,novakNW12,novakPU21,novakTem20}
for related results.
This also implies that they are optimal up to logarithmic factors
for the problem of integration on 
$ H^{k,{{\rm mix}}}(\mathbb T^d)$
with respect to both error criteria.
We do not know whether the logarithmic loss 
can be avoided
with uniformly distributed samples,
as it is the case with isotropic smoothness.
This is another Open Problem. 

Known optimal sampling points for the integration problem
on 
$ H^{k,{{\rm mix}}}(\mathbb T^d)$
have a very particular structure,
see e.g.\ \cite[Section~8.5]{novakDTU16} and the references therein
for the uniform error 
and \cite{novakKN17,novakUl17} for the Monte Carlo error.
\end{remark}

\noindent
{\bf Open Problem: } \ 
In this section we presented results on the optimal 
order of convergence if we may only use iid sample points.  
There are other ways to compare the power of different algorithms.
In particular, one should study ``iid tractability'', 
i.e., 
tractability properties 
of algorithms that are based on iid sample points, 
see also~\cite{novakLW2021,novakNW12}.  
We already know that tractability properties are 
not related to optimal order of convergence, 
so there is much room for new research. 

\medskip

It would be good to know more about the difference between optimal 
and iid sample points, depending on the problem, the dimension $d$ 
and the domain $D_d  \subset {\mathbb{R} }^d$. 
One possibility is to study the asymptotic constants, such as 
\[
 \lim_{n \to \infty}  
 e\left(\mathcal{P}_n,W_p^s(D_d ),L_q(D_d )\right)
\cdot 
 n^{s/d-(1/p-1/q)_+} .
\]
We conjecture that all these asymptotic constants exist and 
possibly they do not depend on the shape of $D_d $,
only on its volume. 
The asymptotic constants are known only in very rare cases, see \cite{novakNo20}, though. 
Here we present an example from \cite{novakHKNPUsurvey} 
that nicely shows the quality of iid samples.
We study $L_1$-approximation for functions from the class 
$$
F_d = \{ f: [0,1]^d \to {\mathbb{R} } \mid 
|f(x)-f(y)| \le d(x,y) \},
$$
with the maximum metric on the $d$-torus, i.e., 
$$
d(x,y) = \min_{k \in \mathbb{Z}^d} \Vert x+k-y \Vert_\infty .
$$
Then
\[
\lim_{n \to \infty} 
{e^{\rm unif}}\left(\mathcal{P}_n,F_d, \novakchange{L_1([0,1]^d)} \right)
\cdot n^{1/d} = 
\frac{1}{2} \, \Gamma (1+ 1/d) \approx  \frac{1}{2} - \frac{\gamma}{2d} 
\]
with the Euler number $\gamma \approx 0.577$, noting that $\Gamma'(1)=-\gamma$. 
This compares very well with the error of optimal methods that, for $n=m^d$,  equals 
$\frac{d}{2d+2} \, n^{-1/d}$. 
To achieve the same error $\varepsilon$ in high dimension 
with random iid sample points, we have to multiply the number 
of optimal sample points by roughly $\exp(1-\gamma)\approx 1.526$;
this factor is quite small and does not 
increase with $d$. 

\section{Concluding Remarks}   

\subsubsection*{Tractability} \ 
We defined the curse of dimensionality and certain notions 
of tractability above, see \eqref{novakcurse}. 
Results till 2012 are summarized in the monographs~\cite{novakNW08,novakNW10,novakNW12} 
and we presented a few recent results. 

Tractability is a very active research area and we mention a few more important papers: 
Stability properties are studied by Dick and Goda~\cite{novakDG20};
the randomized error is studied 
by Dick, Goda and Suzuki~\cite{novakDGS21} and 
by
Kritzer, Kuo, Nuyens and M.~Ullrich~\cite{novakKKNU19};
many (old and new) results can be found in the 
recent book by Dick, Kritzer and Pillichshammer~\cite{novakDKP22},
see also the introduction by 
Leobacher and Pillichshammer~\cite{novakLP14}. 

Dick~\cite{novakDi14}, 
Goda~\cite{novakGo22,novakGo23} 
\novakchange{and Krieg~\cite{novakK23}} found  unweighted function spaces
where polynomial tractability holds;  
the Sobolev space that goes with the star discrepancy is an earlier 
example in this direction, see~\cite{novakHNWW99}. 
\novakchange{The curse of dimensionality for the $L_p$-discrepancy holds 
for $p$ of the form $p=2\ell/(2\ell-1)$ with 
a natural number $\ell$, see~\cite{novakNP23}.}

The Clenshaw Curtis Smolyak algorithm,
introduced and studied  
more than twenty years ago
in~\cite{novakBNR99,novakNR96b,novakNR97,novakNR99} 
has  tractability properties for very smooth functions, 
see~\cite{novakHNU14} 
and Xu~\cite{novakXu15}.  

Traditionally, tractability studies deal with functions of 
finite smoothness and then we can at most hope 
for polynomial tractability.
An important new topic is weighted spaces and exponential tractability
where, usually, analytic functions are studied and one might expect a much 
faster rate of convergence. 
I invite the reader to study the papers
Kritzer and Wo\'zniakowski~\cite{novakKW19},  
Ebert and Pillichshammer~\cite{novakEP21} 
and Krieg, Siedlecki, M. Ullrich and Wo\'zniakowski~\cite{novakKSUW22},
where one  can find more references. 
Also the papers \cite{novakNUWZ18}
and 
Vyb\'\i ral~\cite{novakVy14}
deal with very smooth functions. 

\subsubsection*{Oscillatory integrals} 

Oscillatory integrals are important for applications and are well studied in numerical analysis. 
There exist error bounds but only few complexity results (lower bounds and fitting upper 
bounds). One reason might be that some authors do not fix a class $F$ of inputs and 
then a complexity analysis cannot be done. 
See \cite{novakNUW13,novakNUWZ17,novakZN19}
for complexity results in the univariate case and for further references.  
The paper Wu, Graham, Ma and Zhang~\cite{novakWGMZ22} 
contains interesting results for the multivariate case and 
a Filon-Clenshaw-Curtis-Smolyak rule.

In \cite{novakZN19} 
we study optimal quadrature formulas for 
weighted integrals and integrands from the Sobolev space $H^1([0,1])$.
A particular case is the computation of Fourier coefficients,
where the oscillatory weight
is given by $\rho_k(x) = \exp(- 2 \pi i k x)$.
Here we study the question whether equidistant nodes are optimal
or not. We prove that this depends on $n$ and $k$:
equidistant nodes are \emph{optimal}
if $n \ge  2.7 |k| +1 $ but might be suboptimal for small $n$.
In particular, the equidistant nodes
$x_j = j/ |k|$ for $j=0, 1, \dots , |k|$, hence $n=|k|+1$, 
are the \emph{worst}  possible nodes
and do not give any useful information.

\subsubsection*{A final open problem} 

Study the 
approximation  
of the integral
\begin{equation}\label{novakeq27}  
S_d(f) = \int_{D_d} f(x)  \varrho(x) \, \rd x
\end{equation}
over an open subset $D_d\subset \R^d$  
for integrable functions
$f\colon D_d\to\R$. 
We are interested in universal 
randomized algorithms $(A_n^\omega)$ of the form 
\[
A_n^\omega (f) = \sum_{j=1}^n a_j(\omega) f(x_j(\omega)) \, , 
\]
possibly based on iid uniformly distributed $x_j$. 
Ideally, we would like to have the same optimal error bounds as 
in Theorem~\ref{novakT3}. 
It is not clear whether such algorithms exist, we know it 
only for $\rho=1$ on $D_d=[0,1]^d$. 
In addition,  the optimal order of convergence 
may be useless if the involved constants are huge. 
It might be more practical 
to  relax the order of the error bounds and to  allow 
additional log factors. 
The papers 
\cite{novakKUV21,novakKNS21,novakKU21,novakMN18}
could  be useful to reach this goal. 

One idea would be to use  exactness spaces 
\[
V_K =  {\rm  span} ( \{ x^i \mid \prod_{j=1}^d (i_j +1) \le K \} ) 
\]
for $K=1,2, \dots $ 
and proceed similarly as in \cite{novakMN18}. 

\medskip

I invite the reader to study also my earlier 
survey paper \cite{novakNo16} where one may find more results and references. 

%  \bigskip

%%%%%%%%%%%%%%%%%%%%%%%%%%%%%%%%%%%%%%%%%%%%%%%%%%%%%%%%%%%%%%%%%%%%%%%%

\begin{acknowledgement}
I thank the \novakchange{referee and the} following colleagues and friends for valuable 
remarks: 
David Krieg,  
Peter Kritzer, 
Fritz Pillichshammer, 
Mathias Sonnleitner, 
Mario Ullrich, 
Jan Vyb{\'\i}ral
and Marcin Wnuk.
\end{acknowledgement}

%%%%%%%%%%%%%%%%%%%%%%%%%%%%%%%%%%%%%%%%%%%%%%%%%%%%%%%%%%%%%%%%%%%%%%%
%%% The bibliography
%
% BibTeX users please use
%\bibliographystyle{spmpsci}
%\bibliography{mybibfile}
% and then copy paste the contents of the .bbl file here for the final version.
%
% E.g.:

\end{document}